\documentclass{article}

\usepackage{amsmath}
\usepackage{url}

\usepackage{amsfonts, amsthm, amsmath, amssymb}
\newtheorem{conjecture}{Conjecture}
\newtheorem{case}{Case}
\newtheorem{statement}{Statement}
\theoremstyle{remark}
\newtheorem{rem}{Remark}

\usepackage{hyperref}
\usepackage{breakurl}

\usepackage{mathtools}

\allowdisplaybreaks

\usepackage[nottoc]{tocbibind}

\hypersetup{
pdftitle={On Legendre's, Brocard's, Andrica's, and Oppermann's Conjectures},
pdfauthor={Germ\'an Andr\'es Paz},
pdfsubject={},
pdfkeywords ={Andrica's conjecture, Brocard's conjecture, Legendre's conjecture, Oppermann's conjecture, prime numbers}
}

\begin{document}

\title{On Legendre's, Brocard's, Andrica's, and Oppermann's Conjectures}
\author{Germ\'an Andr\'es Paz}
\date{}

\maketitle

\centerline{Instituto de Educaci\'on Superior N$^\circ$28 Olga Cossettini}
\centerline{(2000) Rosario, Santa Fe, Argentina}
\centerline{E-mail: \texttt{germanpaz\_ar@hotmail.com}}

\centerline{}

\centerline{\bf Abstract}
Let $n\in\mathbb{Z}^+$. Is it true that every sequence of $n$ consecutive integers greater than $n^2$ and smaller than $(n+1)^2$ contains at least one prime number? In this paper we show that this is actually the case for every $n \leq 1,193,806,023$. In addition, we prove that a positive answer to the previous question for all $n$ would imply Legendre's, Brocard's, Andrica's, and Oppermann's conjectures, as well as the assumption that for every $n$ there is always a prime number in the interval $[n,n+2\lfloor\sqrt{n}\rfloor-1]$.



\section{Introduction}\label{Introduction}

The well-known Bertrand's postulate states that for every integer $n>3$ there always exists a prime number $p$ such that $n<p<2n-2$ (another formulation of this theorem is that for every $n>1$ there always exists a prime number $p$ such that $n<p<2n$). This statement, which had been conjectured by Joseph Bertrand in 1845, was first proved by P. L. Chebyshev in 1850. In 1919, Ramanujan \cite{Ramanujan} gave a simpler proof, from which the concept of Ramanujan primes would later arise. Erd\H{o}s published another simple proof of Bertrand's postulate in 1932 \cite{Erdos}. 

After Bertrand's postulate was proved in 1850, better results have been obtained by using both elementary and nonelementary methods. In 1930, Hoheisel \cite{Hoheisel} showed that there exists a constant $\theta<1$ such that 
\begin{equation}\label{asymptotic_formula}
\pi(x+x^\theta)-\pi(x)\sim\frac{x^\theta}{\ln x}\text{,}
\end{equation}
where $\pi$ denotes the prime-counting function. In fact, Hoheisel showed that one may take $\theta=32999/33000$. This result was later improved to $\theta=249/250$ by Heilbronn \cite{Heilbronn}, and to $\theta=3/4+\varepsilon$ by Tchudakoff \cite{Tchudakoff}. In 1937, Ingham \cite{Ingham} proved that there exists a positive integer $K$ such that $p_{n+1}-p_n<K(p_n)^{5/8}$, where $p_n$ is the $n$th prime number. A consequence of Ingham's result is that there is always a prime number between $n^3$ and $(n+1)^3-1$ if $n$ is an integer greater than $K^8$ \cite{Mills}. In 1952, Jitsuro Nagura \cite{Nagura} showed that there exists a prime in the interval $[n,6n/5]$ for every $n\ge 25$. Huxley \cite{Huxley} showed that, for sufficiently large $n$, we have $p_{n+1}-p_n<(p_n)^\delta$ whenever $\delta>7/12$. In 1976, Lowell Schoenfeld \cite{Schoenfeld} proved that the interval $(n,n+n/16597)$ contains a prime for all $n\ge 2,010,760$. Iwaniec and Pintz \cite{Iwaniec_Pintz} showed that there is always a prime in the interval $[n-n^{23/42},n]$ for sufficiently large $n$. In 1998, Pierre Dusart \cite{Dusart} showed in his doctoral thesis that there exists a prime in the interval $[n,n+n/(2 \ln^2 n)]$ for every $n\ge 3,275$. In 2001, Baker, Harman, and Pintz \cite{Baker_Harman_Pintz} proved that in (\ref{asymptotic_formula}) the constant $\theta$ may be taken to be $0.525$. In other words, these authors showed that the interval $[x,x+x^{0.525}]$ contains at least one prime number for sufficiently large $x$. In 2006, M. El Bachraoui \cite{M. El Bachraoui} gave a proof of the fact that there is a prime in the interval $[2n,3n]$. Dusart improved his previous result in 2010 \cite{Dusart2}, when he proved that for $n \ge 396,738$ there exists a prime in the interval $[n,n+n/(25 \ln^2 n)]$. In 2011, Andy Loo \cite{Loo} provided a proof of the fact that there is always a prime number in the interval $[3n,4n]$. Moreover, we can also say that if the Riemann hypothesis is true, then in (\ref{asymptotic_formula}) we can take $\theta=1/2+\varepsilon$ \cite{Maier}.

Although much progress has been made towards finding shorter and shorter intervals containing at least one prime, there are still many open problems in Number Theory regarding the existence of prime numbers in certain intervals. Some of these problems are: 1) Legendre's conjecture (see Sect. \ref{Legendre}), 2) Brocard's conjecture (Sect. \ref{Brocard}), 3) Andrica's conjecture (Sect. \ref{Andrica}), and 4) Oppermann's conjecture (Sect. \ref{Oppermann}).

In this paper we consider the following conjecture regarding the distribution of prime numbers:

\begin{conjecture}\label{Conjecture1}
If $n$ is any positive integer and we take $n$ consecutive integers
located between ${n^2}$ and ${\left( {n + 1} \right)^2}$, then among those
$n$ integers there is at least one prime number. In other words, if $a_{1}$, $a_{2}$, $a_{3}$, $a_{4}$, \ldots{}, $a_{n}$ are $n$ consecutive integers
such that $n^2<a_1<a_2<a_3<a_4< ... <a_n<(n + 1)^2$, then at least one of those $n$ integers is a prime number.
\end{conjecture}

\begin{rem}
Throughout this paper, whenever we say that a number $b$ is
\emph{between} a number $a$ and a number $c$, it means that
$a < b < c$, which means that $b$ is never equal to $a$ or
$c$. Moreover, the number
$n$ that we use in this document is always a positive
\emph{integer}.
\end{rem}
While it is conjectured that (\ref{asymptotic_formula}) holds for all $\theta\in(0,1)$, it would be interesting to find a relation that explains why a sequence of $n$ consecutive integers greater than $n^2$ and smaller than $(n+1)^2$ cannot contain only composite numbers. In other words, it would be interesting to find a relation between the amount of consecutive integers in a sequence and the perfect squares between which the mentioned sequence is located.

Let us see some cases in which Conjecture \ref{Conjecture1} is true:

\begin{itemize}

\item If we consider $n=1$, we have $1^2<2<3<(1+1)^2$, and we can see that the numbers 2 and 3 are both prime numbers.

\item If we consider $n=2$, we have $2^2<5<6<7<8<(2+1)^2$. If we take any sequence of 2 consecutive integers greater than $2^2$ and smaller than $(2+1)^2$, then at least one of those 2 integers is a prime number. This is true because each of the sequences $\{5, 6\}$, $\{6, 7\}$, and $\{7, 8\}$ contains at least one prime number.

\item If we consider the case where $n=3$, we have $3^2<10<11<12<13<14<15<(3+1)^2$. It is easy to verify that each of the sequences $\{10, 11, 12\}$, $\{11, 12, 13\}$, $\{12, 13, 14\}$, and $\{13, 14, 15\}$ contains at least one prime number.

\end{itemize}
We can easily prove that Conjecture \ref{Conjecture1} is also true for $n=4$, $n=5$, $n=6$, and larger values of $n$. 

Let us suppose that $p$ and $q$ are two consecutive prime numbers such that $p<q$. It is easy to verify that the amount of composite numbers between $p$ and $q$ is equal to $q-p-1$. In other words, the gap between $p$ and $q$ is equal to $q-p-1$.

A gap between two consecutive prime numbers is said to be \emph{maximal} if it is larger than all gaps between smaller primes. Let $P$ denote a prime number followed by a maximal gap $G$. This means that $G$ is a maximal gap between two consecutive primes $P$ and $Q$. Taking into account that $G=Q-P-1$, we can use the following method to find out the values of $n$ for which Conjecture \ref{Conjecture1} is true:

\begin{itemize}

\item We calculate $\lfloor\sqrt{P}\rfloor$.

\item If $G<\lfloor\sqrt{P}\rfloor$, then Conjecture \ref{Conjecture1} holds for every $n$ such that $G-1\leq n \leq\lfloor\sqrt{P}\rfloor-1$ (in this paper, the gap between two consecutive prime numbers greater than 2 is an odd number).

\item If $G=\lfloor\sqrt{P}\rfloor$ and $(\lfloor\sqrt{P}\rfloor+1)^2-P-1<\lfloor\sqrt{P}\rfloor$, then Conjecture \ref{Conjecture1} holds for $n=\lfloor\sqrt{P}\rfloor-1$ and for $n=\lfloor\sqrt{P}\rfloor$.

\item If $G>\lfloor\sqrt{P}\rfloor$, $(\lfloor\sqrt{P}\rfloor+1)^2-P-1<\lfloor\sqrt{P}\rfloor$, and the previous maximal gap is less than $\lfloor\sqrt{P}\rfloor$, then Conjecture \ref{Conjecture1} holds for $n=\lfloor\sqrt{P}\rfloor$.

\end{itemize}
By using this method in combination with tables of maximal gaps (see \cite{Andersen} and \cite{Nicely}) and lists of prime numbers, we verify that Conjecture \ref{Conjecture1} holds at least for every $n$ such that $1\leq n\leq 1,193,806,023$ (we have $10^9<1,193,806,023$). In other words, the conjecture holds over the first $1,193,806,024^2=1,425,172,822,938,688,576$ positive integers. This number is greater than $10^{18}$.

\begin{rem}
Andersen \cite{Andersen} and Nicely \cite{Nicely} define gaps between consecutive primes $p$ and $q$ as $g=q-p$. In this paper, we define a gap between consecutive primes as $g=q-p-1$. This means that in order to work with the method described before in combination with Andersen and Nicely's tables of maximal gaps, we need to substract 1 from the gaps listed on those tables. Therefore, a maximal gap of $x$ in Andersen and Nicely's tables is considered as a maximal gap of $x-1$ in our paper.
\end{rem}

In Sects. \ref{Legendre}, \ref{Brocard}, \ref{Andrica}, and \ref{Oppermann} we will show that if Conjecture \ref{Conjecture1} is true, then Legendre's, Brocard's, Andrica's, and Oppermann's conjectures follow. Moreover, we will also show that if the mentioned conjecture holds, then there exists a prime in the interval $\left[{n,n + 2\left\lfloor{\sqrt n}\right\rfloor - 1}\right]$ for every positive integer $n$ (see Sect. \ref{Interval}).

\section{Legendre's Conjecture}\label{Legendre}

Legendre's conjecture \cite{Legendre} states that for every positive integer $n$ there exists at least one prime number $p$ such that $n^2<p<(n+1)^2$. This conjecture is one of Landau's problems \cite{Hardy_Wright, Landau}.

It is easy to verify that the amount of integers located between ${n^2}$
and ${\left( {n + 1} \right)^2}$ is equal to $2n$.

\begin{proof}
We have
\begin{align*}
{\left( {n + 1} \right)^2} - {n^2} &= 2n + 1\\
{n^2} + 2n + 1 - {n^2} &= 2n + 1\\
2n + 1 &= 2n + 1\text{.}
\intertext{We need to exclude the number ${\left( {n + 1} \right)^2}$ because we are taking into consideration the integers that are greater than ${n^2}$ and smaller than
$\left( {n + 1} \right)^2$. Therefore, we get}
2n + 1 - 1 &= 2n\text{.} \qedhere
\end{align*}
\end{proof}

According to this result, between ${n^2}$ and ${\left( {n + 1} \right)^2}$ there are
two groups of $n$ consecutive integers each that do not have any integer
in common. Example for $n = 3$:
\[
{3^2}\, \underbrace {\underbrace {10\;\quad11\;\quad
12}_{\scriptstyle{\rm{\hspace{25pt}Group\;A }}\hfill\atop
\scriptstyle{\rm{(}}n\;{\rm{consecutive\;integers)}}\hfill}\underbrace
{13\;\quad14\;\quad15}_{\scriptstyle{\rm{\hspace{25pt}Group\;B }}\hfill\atop
\scriptstyle{\rm{(}}n\;{\rm{consecutive\;integers)}}\hfill}}_{2n\;{\rm{consecutive\;
integers}}} {(3 + 1)^2}
\]
\text{Group A} and \text{Group B} do not have any integer in common. Now, according to Conjecture \ref{Conjecture1}, Group A contains at least one prime
number and \mbox{Group B} also contains at least one prime number, which implies that
between ${3^2}$ and ${\left( {3 + 1} \right)^2}$ there are at least \text{two}
prime numbers. This is true because the numbers 11 and 13 are both prime.

All this means that if Conjecture \ref{Conjecture1} is true, then there are at least \text{two} prime numbers between ${n^2}$ and ${\left( {n + 1} \right)^2}$ for every positive
integer $n$. \text{As a result, if Conjecture \ref{Conjecture1} is true, then Legendre's
conjecture is also true.}

\section{Brocard's Conjecture}\label{Brocard}

Brocard's conjecture \cite{Brocard} states that if ${p_n}$ and ${p_{n + 1}}$ are consecutive prime numbers greater than 2, then between ${\left( {{p_n}} \right)^2}$ and ${\left( {{p_{n + 1}}} \right)^2}$ there are at least four prime numbers.

Since $2 < {p_n} < {p_{n + 1}}$, we have ${p_{n + 1}} - {p_n} \ge 2$. This means
that there is at least one positive integer $a$ such that ${p_n} < a < {p_{n +
1}}$. As a result, there exists at least one positive integer \textit{a }such that ${\left( {{p_n}} \right)^2} < {a^2} < {\left( {{p_{n + 1}}}
\right)^2}$.

Conjecture \ref{Conjecture1} states that between ${\left( {{p_n}} \right)^2}$ and ${a^2}$ there are at least two prime numbers and that between ${a^2}$ and ${\left( {{p_{n +
1}}} \right)^2}$ there are also at least two prime numbers. In other words, if
Conjecture \ref{Conjecture1} is true, then there are at least four prime numbers between ${\left(
{{p_n}} \right)^2}$ and ${\left( {{p_{n + 1}}} \right)^2}$. As a consequence, if Conjecture \ref{Conjecture1} is true, then Brocard's conjecture is also true.

\section{Andrica's Conjecture}\label{Andrica}

Andrica's conjecture \cite{Andrica1, Andrica2} states that $\sqrt {{p_{n + 1}}}  - \sqrt {{p_n}}  < 1$ for
every pair of consecutive prime numbers ${p_n}$ and ${p_{n + 1}}$ (of course, ${p_n} < {p_{n + 1}}$).

Obviously, every prime number is located between two consecutive perfect
squares. Now, let us suppose that $p$ is any prime number and that $q$ is the prime number immediately following $p$. If we take into account that $p$ is obviously located between
${n^2}$ and ${\left( {n + 1} \right)^2}$ for some $n$, two things may happen:

\begin{case}\label{case1}
The number $p$ is among the first $n$ consecutive
integers that are located between ${n^2}$ and ${\left( {n + 1} \right)^2}$. These
$n$ integers form what we call `Group A,' and the following
$n$ integers form what we call `Group B.'
\end{case}
Let us look at the following graphic.
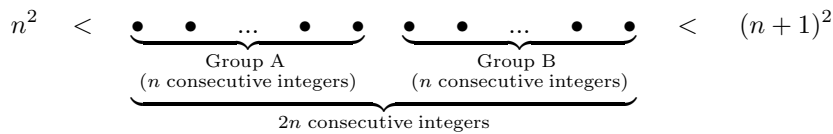
\[
{n^2}\quad  < \,\quad \underbrace {\underbrace { \bullet \;\quad  \bullet
\;\quad ...\;\quad  \bullet \;\quad  \bullet }_{\scriptstyle{\rm{\hspace{25pt}Group\;A}}\hfill\atop
\scriptstyle{\rm{(}}n{\rm{\;consecutive\;integers)}}\hfill}\;\quad \underbrace {
\bullet \;\quad  \bullet \;\quad ...\;\quad  \bullet \;\quad  \bullet
}_{\scriptstyle{\rm{\hspace{25pt}Group\;B}}\hfill\atop
\scriptstyle{\rm{(}}n{\rm{\;consecutive\;integers)}}\hfill}}_{2n\;{\rm{consecutive\;
integers}}}\;\quad  < \;\quad {(n + 1)^2}
\]
If $p$ is located in Group A and Conjecture \ref{Conjecture1} is true, then $q$ is
either located in Group A or in Group B. In both cases we have $\sqrt{q}-\sqrt{p}<1$, since $\sqrt {{{(n + 1)}^2}}  - \sqrt {{n^2}}  = 1$ and the numbers $\sqrt
{q}$ and $\sqrt {p}$ are closer to each other than $\sqrt {{{(n +
1)}^2}} $ in relation to $\sqrt {{n^2}} $.

\begin{case}\label{case2}
The prime number $p$ is located in Group B.
\end{case}

If $p$ is located in Group B and Conjecture \ref{Conjecture1} is true, it may happen that
$q$ is also located in Group B. In this case, it is very easy to verify
that $\sqrt {q}  - \sqrt {p}  < 1$, as explained before.

Otherwise, if $q$ is not located in Group B, then $q$ is located in `Group C.' In this case, the largest value that $q$ can have is \mbox{$q={(n + 1)^2} + n + 1={n^2} +
2n + 1 + n + 1={n^2} + 3n + 2$}, while the smallest value $p$ can have is $p = {n^2} + n + 1$ (in order to make the process easier, we are
not taking into account the fact that in this case the numbers $p$ and $q$
have different parity, which means that they cannot be both prime at the same time).

This means that the largest possible difference between $\sqrt {q} $
and $\sqrt {p} $ is $\sqrt {q}  - \sqrt {p}  = \sqrt {{n^2} +
3n + 2}  - \sqrt {{n^2} + n + 1} $.
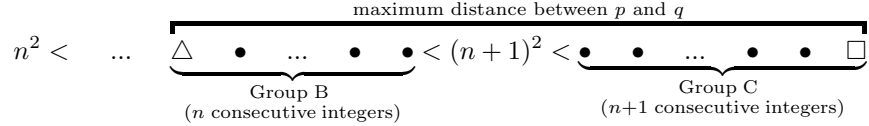
\[
{n^2}< \;\quad ... \;\quad \overbracket[1.5pt]{\underbrace { \triangle \;\quad  \bullet
\;\quad ...\;\quad  \bullet \;\quad  \bullet }_{\scriptstyle{\rm{\hspace{25pt}Group\;B}}\hfill\atop
\scriptstyle{\rm{(}}n{\rm{\;consecutive\;integers)}}\hfill}<{(n + 1)^2}<\underbrace { \bullet
\;\quad  \bullet \;\quad ...\;\quad  \bullet \;\quad  \bullet \;\quad  \Box
}_{\scriptstyle{\rm{\hspace{25pt}\;\;Group\;C}}\hfill\atop
\scriptstyle{\rm{(}}n + 1{\rm{\;consecutive\;integers)}}\hfill}}^\text{maximum distance between $p$ and $q$}
\]

\raggedright{
$\triangle={n^2} + n + 1 = p$\\
$\Box={n^2} + 3n + 2 = q$
}

\hspace{15.0pt}It is easy to prove that $\sqrt {{n^2} + 3n + 2}  - \sqrt {{n^2} + n + 1}  < 1$.

\begin{proof}
We have
\begin{align*}
\sqrt {{n^2} + 3n + 2}  - \sqrt {{n^2} + n + 1}   &<  1\\
\sqrt {{n^2} + 3n + 2}   &<  1 + \sqrt {{n^2} + n + 1}\\
{n^2} + 3n + 2  &<  {\left( {1 + \sqrt {{n^2} + n + 1} } \right)^2}\\
{n^2} + 3n + 2  &<  1 + 2\sqrt {{n^2} + n + 1}  + {n^2} + n + 1\\
{n^2} + 3n + 2 - {n^2} - n - 1  &<  1 + 2\sqrt {{n^2} + n + 1}\\
2n + 1  &<  1 + 2\sqrt {{n^2} + n + 1}\\
2n  &<  2\sqrt {{n^2} + n + 1}\\
n  &<  \frac{{2\sqrt {{n^2} + n + 1} }}{2}\\
n  &<  \sqrt {{n^2} + n + 1}\\
{n^2}  &<  {n^2} + n + 1,
\end{align*}
which is true for every positive integer $n$. \qedhere
\end{proof}

\parbox[c]{345.0pt}{
\begin{rem}
In general, to prove that an inequality is correct, we can solve that inequality step by step. If we get a result which is obviously correct, then we can start with that correct result, `work backwards from there' and prove that the initial statement is true.
\end{rem}}

\parbox[c]{345.0pt}{\hspace{15.0pt}We can see that even when the difference between $q$ and $p$ is the largest possible difference, we have $\sqrt{q}-\sqrt{p}<1$. If the difference between $q$ and $p$ were smaller, then of course it would also happen that $\sqrt{q}-\sqrt{p}<1$.}

\parbox[c]{345.0pt}{\hspace{15.0pt}According to Cases \ref{case1} and \ref{case2}, if Conjecture \ref{Conjecture1} is true, then Andrica's conjecture is also true.}

\section{Oppermann's Conjecture}\label{Oppermann}

\parbox[c]{345.0pt}{
Oppermann's conjecture \cite{Oppermann} states that for any integer $n>1$ there is a prime number in the interval $[n^2-n,n^2]$ and another prime in the interval $[n^2,n^2+n]$. Now, if Conjecture \ref{Conjecture1} holds, then there exists a prime in the interval $[n^2-n+1,n^2-1]$ and another prime in the interval $[n^2+1,n^2+n]$, for $n>1$. This means that if Conjecture \ref{Conjecture1} is true, then so is Oppermann's conjecture.
}

\vspace{2pt}
\parbox[c]{345.0pt}{
\hspace{15.0pt}Note that according to Oppermann's conjecture there is a prime in the interval $[n^2-n+1,n^2-1]$ and another prime in the interval $[n^2+1,n^2+n]$ for $n>1$, whereas we ask the question whether \emph{any} sequence of $n$ consecutive integers greater than $n^2$ and smaller than $(n+1)^2$ contains at least one prime number. In other words, if Conjecture \ref{Conjecture1} is true, then Oppermann's conjecture follows, whereas the reciprocal is not necessarily true.
}

\section{On the Interval $\left[n,n+2\left\lfloor\sqrt{n}\right\rfloor-1\right]$}\label{Interval}

\parbox[c]{345.0pt}{It is easy to verify that if Conjecture \ref{Conjecture1} is true, then in the interval $\left[
{{n^2} + n + 1,{n^2} + 3n + 2} \right]$ there are at least two prime numbers for every positive integer $n$.}

\hspace{15.0pt}Now, the number ${n^2} + n + 1$ is always an odd integer.

\begin{proof}
\[
\]
\begin{itemize}
\item If $n$ is even, then ${n^2}$ is also even. Then we have
\[
(even\;integer + even\;integer) + 1 = even\;integer +
odd\;integer = 
\]
\[
= odd\;integer.
\]
\item If $n$ is odd, then ${n^2}$ is also odd. Then we have
\[
(odd\;integer + odd\;integer) + 1 = even\;integer +
odd\;integer = 
\]
\[
= odd\;integer. \qedhere
\]
\end{itemize}
\end{proof}

\parbox[c]{345.0pt}{\hspace{15.0pt}Since the number ${n^2} + n + 1$ is always an odd integer, then it may be
prime or not. Now, the number ${n^2} + 3n + 2$ can never be prime, since this number is always an even integer greater than 2.}

\begin{proof}
\[
\]
\begin{itemize}
\item If $n = 1$ (smallest value $n$ can have), then ${n^2} + 3n + 2 = 1 + 3 +
2 = 6$.
\item If $n$ is even, then ${n^2}$ and $3n$ are both even integers. The number
\mbox{2} is also an even integer, and we know that
\[
even\;integer + even\;integer + even\;integer = even\;integer.
\]
\item If $n$ is odd, then ${n^2}$ and $3n$ are both odd integers, and we know
that
\[
(odd\;integer + odd\;integer) + even\;integer = 
\]
\[
= even\;integer +
even\;integer = even\;integer. \qedhere
\]
\end{itemize}
\end{proof}

\parbox[c]{345.0pt}{\hspace{15.0pt}From all this we deduce that if Conjecture \ref{Conjecture1} is true, then the maximum distance
between two consecutive prime numbers is the one from the number ${n^2} + n + 1$
to the number ${n^2} + 3n + 2 - 1 = {n^2} + 3n + 1$, which means that
in the interval $[{n^2} + n + 1,{n^2} + 3n + 1]$ there are at least two prime
numbers. In other words, in the interval $[{n^2} + n + 1,{n^2} + 3n]$ there
is at least one prime number.}
\parbox[c]{345.0pt}{\hspace{15.0pt}The difference between the numbers ${n^2} + n + 1$ and ${n^2} + 3n$ is \mbox{${n^2} +
3n - ({n^2} + n + 1) = {n^2} + 3n - {n^2} - n - 1 = 2n - 1$}. In addition to this,
$\left\lfloor {\sqrt {{n^2} + n + 1} } \right\rfloor  = n$. This means that in the interval $\left[ {{n^2} + n + 1,{n^2} + n + 1 +
2\left\lfloor {\sqrt {{n^2} + n + 1} } \right\rfloor  - 1} \right]$ there is at least one prime number. In other words, if $a = {n^2} + n + 1$, then the
interval $\left[ {a,a + 2\left\lfloor {\sqrt a } \right\rfloor  - 1} \right]$
contains at least one prime number.}
\parbox[c]{345.0pt}{
\begin{rem}
The symbol $\left\lfloor {} \right\rfloor $ represents the \textit{floor function}. The floor function of a given number is the largest integer that is not greater than that number. For example, $\left\lfloor 3.5 \right\rfloor=3$.
\end{rem}
}
\parbox[c]{345.0pt}{\hspace{15.0pt}Now, if Conjecture \ref{Conjecture1} is true, then the following statements are all true:}
\parbox[c]{345.0pt}{
\begin{statement}\label{statement1}
If $a$ is a perfect square, then in the interval $\left[ {a,a +
\left\lfloor {\sqrt a } \right\rfloor } \right]$ there is at least one prime
number.
\end{statement}}
\parbox[c]{345.0pt}{
\begin{statement}\label{statement2}
If $a$ is an integer such that ${n^2} < a \le {n^2} + n + 1 < {(n +
1)^2}$, then in the interval $\left[ {a,a + \left\lfloor {\sqrt a }
\right\rfloor  - 1} \right]$ there is at least one prime number.
\end{statement}}
\parbox[c]{345.0pt}{
\begin{statement}\label{statement3}
If $a$ is an integer such that ${n^2} < {n^2} + n + 2 \le a < {(n +
1)^2}$, then in the interval $\left[ {a,a + 2\left\lfloor {\sqrt a }
\right\rfloor  - 1} \right]$ there is at least one prime number.
\end{statement}}

\parbox[c]{345.0pt}{\hspace{15.0pt}We know that $a + 2\left\lfloor {\sqrt a } \right\rfloor  - 1 \ge a +
\left\lfloor {\sqrt a } \right\rfloor $.}

\begin{proof}
We have
\begin{align*}
a + 2\left\lfloor {\sqrt a } \right\rfloor  - 1 &\ge a + \left\lfloor {\sqrt a }
\right\rfloor\\   
2\left\lfloor {\sqrt a } \right\rfloor
 - 1 &\ge \left\lfloor {\sqrt a } \right\rfloor\\
2\left\lfloor {\sqrt a } \right\rfloor  &\ge \left\lfloor {\sqrt a } \right\rfloor
 + 1\\
\left\lfloor {\sqrt a } \right\rfloor  + \left\lfloor
{\sqrt a } \right\rfloor  &\ge \left\lfloor {\sqrt a } \right\rfloor  + 1\\
\left\lfloor {\sqrt a } \right\rfloor  &\ge 1\text{,}
\end{align*}
which is true for every positive integer $a$. \qedhere
\end{proof}

\hspace{15.0pt}And we also know that $a + 2\left\lfloor {\sqrt a } \right\rfloor  - 1 > a +
\left\lfloor {\sqrt a } \right\rfloor  - 1$.

\begin{proof}
We have
\begin{align*}
a + 2\left\lfloor {\sqrt a } \right\rfloor  - 1 &> a + \left\lfloor {\sqrt a }
\right\rfloor  - 1\\
2\left\lfloor {\sqrt a }
\right\rfloor  &> \left\lfloor {\sqrt a } \right\rfloor\text{,}
\end{align*}
which is obviously true for every positive integer $a$. \qedhere
\end{proof}

\parbox[c]{345.0pt}{\hspace{15.0pt}All this means that the interval $\left[ {a,a + 2\left\lfloor {\sqrt a }
\right\rfloor  - 1} \right]$ can be applied to the number $a$ in Statement \ref{statement1}, to the number $a$ in Statement \ref{statement2}, and to the number
$a$ in Statement \ref{statement3}. Therefore, if $n$ is \text{any} positive integer and Conjecture \ref{Conjecture1} is true, then in the interval $\left[ {n,n + 2\left\lfloor {\sqrt n } \right\rfloor  - 1}
\right]$ there is at least one prime number (in order to provide more standardized notation, we are now replacing letter {$a$} with letter {$n$}). According to this, we can also say that if Conjecture \ref{Conjecture1} is
true, then there is always a prime number in the interval $\left[ {n,n + 2\sqrt n  - 1} \right]$ for every positive integer $n$.
}

\section{Conclusion}

\vspace{4pt}

\parbox[c]{345.0pt}{We have proved that if Conjecture \ref{Conjecture1} is true, then Legendre's, Brocard's, Andrica's, and Oppermann's conjectures follow. In addition, we have shown that if the mentioned conjecture holds, then the interval $\left[{n,n + 2\left\lfloor{\sqrt n}\right\rfloor - 1}\right]$ contains a prime for every positive integer $n$.} 

\parbox[c]{345.0pt}{\hspace{15.0pt}Now, the number $2\lfloor\sqrt{n}\rfloor-1$ is always an odd number. Suppose that Conjecture \ref{Conjecture1} is true and that $b$ is an odd integer greater than 1. Then the interval $[{b,b + 2\lfloor{\sqrt b}\rfloor - 2}]$ contains at least one prime number. Under the same assumption we can state that if $c$ is any positive even integer, then there exists a prime in the interval $[{c+1,c + 2\lfloor{\sqrt c}\rfloor - 1}]$.}

\parbox[c]{345.0pt}{\hspace{15.0pt}Conjecture \ref{Conjecture1} implies that if $n$ is any positive integer greater than 1, then at least one of the intervals $[n-\lfloor\sqrt{n}\rfloor+1,n]$ and $[n,n+\lfloor\sqrt{n}\rfloor-1]$ contains a prime. Note that we can state the following:

\begin{itemize}

\item If $d$ is a certain positive odd integer such that $d$ is not a perfect square and $\lfloor\sqrt{d}\rfloor$ is even, then the fact that there is a prime in $[d-\sqrt{d},d]$ does not imply that there is a prime in $[d-\lfloor\sqrt{d}\rfloor+1,d]$. Similarly, the fact that there exists a prime in $[d,d+\sqrt{d}]$ does not imply that there exists a prime in $[d,d+\lfloor\sqrt{d}\rfloor-1]$. 

\item If $f$ is a certain positive even integer such that $f$ is not a perfect square and $\lfloor\sqrt{f}\rfloor$ is odd, then the fact that there is a prime in $[f-\sqrt{f},f]$ does not imply that there is a prime in $[f-\lfloor\sqrt{f}\rfloor+1,f]$. Similarly, the fact that there exists a prime in $[f,f+\sqrt{f}]$ does not imply that there exists a prime in $[f,f+\lfloor\sqrt{f}\rfloor-1]$.

\end{itemize}
This means that if $n$ is not a perfect square and $n$ and $\lfloor\sqrt{n}\rfloor$ have different parity, then the fact that $\pi[n-\sqrt{n},n]\ge 1$ does not imply that $\pi[n-\lfloor\sqrt{n}\rfloor+1,n]\ge 1$, and the fact that $\pi[n,n+\sqrt{n}]\ge 1$ does not imply that $\pi[n,n+\lfloor\sqrt{n}\rfloor-1]\ge 1$.
}


%
%
%

\end{document}